\documentclass{amsart}
\newtheorem{theorem}{Theorem}[section]
\newtheorem{lemma}[theorem]{Lemma}

\newtheorem{definition}[theorem]{Definition}
\newtheorem{example}[theorem]{Example}

\newtheorem{corollary}[theorem]{Corollary}
\newtheorem{proposition}[theorem]{Proposition}

\newtheorem{remark}[theorem]{Remark}

\numberwithin{equation}{section}

\def\C{\mathbb C}
\def\R{\mathbb R}
\def\X{\mathbb X}

\def\Y{\mathbb Y}

\def\N{\mathbb N}
\def\cal{\mathcal}
\def\cD{\cal D}

\def\F{\cal F}

\begin{document}
\title[Circular spectrum and bounded solutions]{Circular spectrum and bounded solutions of periodic evolution equations}
\author{Nguyen Van Minh}
\address{Department of Mathematics, University of West Georgia, Carrollton, GA 30118, USA}
\email{vnguyen@westga.edu}
\author{Gaston N'guerekata}
\address{Department of Mathematics, Morgan State University, 1700 E. Cold Spring Lane,
Baltimore, MD 21251, USA} \email{gnguerek@jewel.morgan.edu}
\author{Stefan Siegmund}
\address{Department of Mathematics, Dresden University of Technology, 01062 Dresden, Germany}
\email{siegmund@tu-dresden.de}
\thanks{The authors thank the anonymous referee for his carefully reading the manuscript, pointing out some inacurracies and making several
suggestions to improve
the presentation of the paper.}
\date{\today}
\begin{abstract}We consider the existence and uniqueness of bounded
solutions of periodic evolution equations of the form
$u'=A(t)u+\epsilon H(t,u)+f(t)$, where $A(t)$ is, in general, an
unbounded operator depending $1$-periodically on $t$, $H$ is
1-periodic in $t$, $\epsilon$ is small,
and $f$ is a bounded and continuous function that is not necessarily
uniformly continuous. We propose a new approach to the spectral
theory of functions via the concept of "circular spectrum" and then
apply it to study the linear equations $u'=A(t)u+f(t)$ with
general conditions on $f$. For small $\epsilon$ we show that the
perturbed equation inherits some properties of the linear
unperturbed one. The main results extend recent results in the direction,
saying that if the unitary spectrum of the monodromy operator does
not intersect the circular spectrum of $f$, then the evolution
equation has a unique mild solution with its circular spectrum
contained in the circular spectrum of $f$.
\end{abstract}

\maketitle


\section{Introduction}
The main concern of this paper is the existence and uniqueness of
bounded solutions to periodic evolution equations of the form
\begin{equation}\label{peve}
\frac{du }{dt}= A(t)u  +f(t), \quad t\in\R ,
\end{equation}
and nonlinear perturbed equations of the form
\begin{equation}\label{perturb peve}
\frac{du }{dt}= A(t)u +\epsilon H(t,u) +f(t), \quad t\in\R ,
\end{equation}
where $A(t)$ is, in general, an  unbounded linear operator on a Banach
space $\X$, depending periodically on $t$, $H$ is periodic in $t$
with the same period as $A$, $\epsilon$ is small, and $f$ is in
$L^\infty (\R,\X)$. This is a central problem of the theory and
applications of differential equations. The reader is referred to
\cite{arebathieneu,bardeagen,genbardea,hinnaiminshi,levzhi} and
their references for more information.

\bigskip
Eq. (\ref{perturb peve}) may serve as models for the following
equations
\begin{equation}
\ddot{x}+\omega (t) x + \epsilon h(t,x,\dot x)=f(t),
\end{equation}
where $\omega (t)$ is an $1$-periodic continuous real function,
$h(t,x,\dot x)$ is real continuous, $1$-periodic in $t$ and
uniformly Lipschitz in $(x,\dot x)$ such that $h(t,0,0)=0$, and the
forcing term $f(t)$ is almost periodic, or bounded.

\medskip
It also serves as an abstract setting for the following partial
differential equations (the reader is referred to
\cite{arebathieneu,hen,paz} for more details):
\begin{equation}\label{exa 2}
\begin{cases}  w_t(x,t)=w_{xx}(x,t)+a(t)w(x,t) + \epsilon
 b(t)w(x,t)^2+f(x,t), \ \ \cr \hspace{5cm}0\le s \le \pi ,\ t\ge 0, \cr
w(0,t)=w(\pi ,t)=0 , \ \forall t > 0,
\end{cases}
\end{equation}
where $a(t), b(t),w(x,t), f(x,t)$ are scalar-valued functions,
$a(t),b(t)$ are $1$-periodic and continuous in $t$, and $f(\cdot
,t)$ as an element in $L^2[0,1]$ is almost automorphic. We define
the space ${\X}:=L^2[0,\pi ]$ and $A_T: D(A_T)\subset {\X}\to {\X}$
by the formula
\begin{equation}\label{exa 3}
\begin{cases}\hspace{.5cm} A_Ty=y'', \cr D(A_T)=\{ y\in {\X}:\
y, y' \mbox{ are absolutely continuous}, \ y''\in {\X},\cr
\hspace{4cm} y(0)=y(\pi )=0\}.
\end{cases}
\end{equation}
We define $A(t):= A_T +a(t)$, where $a(t)$ is the operator of
multiplication $a(t)v(\cdot )$ in $\X$, and $H(t,v)(\cdot ):=
b(t)v^2(\cdot )$ for each $v\in\X$. The evolution equation we are
concerned with in this case is the following
\begin{equation}
{\frac{du(t)}{dt}}=A(t)u(t)+\epsilon H(t,u(t)) + f(t), \ u(t) \in
{\X}.
\end{equation}
As is well known (\cite{paz}), the linear part associated with this
evolution equation is well posed, that is, $A(t)$ generates an
$1$-periodic evolutionary process that is strongly continuous.
Therefore, one may include this equation into our abstract model
(\ref{perturb peve}).

\bigskip
In \cite{vu} a conjecture on the existence and uniqueness of an
almost periodic mild solution to (\ref{peve}) is stated when the
unitary spectrum of the  Poincare map $P$ associated with
(\ref{peve}) does not intersect the set $\overline{e^{isp (f)}}$,
where $sp(f)$ denotes the spectrum of $f$ (see the definition
below), and the overline means the closure in the topology of the complex plane. The evolution semigroup method, proposed in \cite{naimin}
(see also \cite{bathutrab,murnaimin,hinnaiminshi}) gives rise to a
positive answer to the conjecture. To use the semigroup theory
machinery a crucial requirement on the strong continuity of the
associated evolution semigroup is made. For instance, $f$ is a
bounded, uniformly continuous function with pre-compact range. In
particular, if $f$ is almost periodic, the condition is
automatically satisfied. In our more general setting $f$ is merely a
bounded and continuous function, so the strong continuity of the
associated evolution semigroup is actually not assumed.
Consequently, the beautiful results of Semigroup Theory do not
apply. This general setting of the problem seems to be natural when
one considers $f$ from some frequently met classes of functions such
as almost automorphic functions (see \cite{bersie,ngu,yi}). On the
other hand, the above-mentioned requirement on $f$ appears to be
technical, and is an obstacle for potential applications of the
results to other areas such as Control Theory. For complete accounts
of results concerned with periodic evolution equations we refer the
reader to
\cite{bathutrab,danmed,hinnaiminshi,levzhi,minnaingu,shinai}. The
asymptotic behavior of evolution equations and applications can be
found in the monographs
\cite{arebathieneu,danmed,hen,levzhi,nee,pru}. For more information
on the spectral theory of functions and its applications the reader
is referred to the monographs
\cite{arebathieneu,kat,levzhi,hinnaiminshi,pru} as well as papers
\cite{imm,min,min2,vusch} and their references.

\bigskip
In this paper we propose a new approach that is based on a concept
of the so-called {\it circular spectrum} of a bounded function $g$
(denoted by $\sigma (g)$). In turn, the circular spectrum of a
bounded function $g$ is defined through a new transform of $g$,
namely, $R(\lambda ,S)g$, where $S$ is the translation to the period
of the coefficient operator $A(\cdot )$ (that is assumed to be $1$).
When the function $g$ is bounded and uniformly continuous, the Weak
Spectral Mapping Theorem in Semigroup Theory yields that $\sigma (g)
= \overline{e^{isp(g)}}$, where $sp(g)$ is the Carleman spectrum of
$g$. This makes our results new extensions of the previous ones to
the general setting. Moreover, a perturbation theory of the linear
equations is given.

\bigskip
Before concluding this introduction section we give an outline of
the paper. We briefly list main notations in Section 2. This section
also contains the definitions as well as properties of almost
periodic and almost automorphic functions. Section 3 deals with a
similar problem for difference equations with continuous time.
Section 4 contains the main result of the paper that deals with the
existence and uniqueness of bounded mild solutions of periodic
evolution equations.

\section{Preliminaries}
\subsection{Basic Notations}
In the paper $\X$ denotes a complex Banach space. The space of all
bounded linear operators in $\X$ is denoted by $L(\X)$ with usual
operator norm. If $A$ is a linear operator (not necessarily bounded)
acting on $\X$, $\sigma (A)$ ($\rho (A)$, respectively) denotes its
spectrum (resolvent set, respectively). The part of spectrum of an
operator $B\in L(\X)$ on the unit circle is denoted by
$\sigma_\Gamma (B)$. For $\lambda\in \rho (A)$, $R(\lambda ,A)$
denotes $(\lambda -A)^{-1}$. If  a mapping $T$ from a Banach space
$\X$ to another Banach space $\Y$ is Lipschitz continuous, then
$$
Lip (T):= \inf \{ L: \| Tx-Ty\| \le L\| x-y\|, \ \mbox{for all}\ x,
y \in \X \}.
$$

In this paper we also use the following notations:
\begin{enumerate}
\item If $z$ is a complex number, then $\Re z$, $\Im z$ denote its
real part and imaginary part, respectively;
\item $BC(\R,\X)$ is the space of all $\X$-valued bounded and
continuous functions on $\R$; $BUC(\R,\X)$ is the function space of
all $\X$-valued bounded and uniformly continuous functions on $\R$;
$KBUC(\R,\X)$ is the function space of all $\X$-valued bounded and
uniformly continuous functions on $\R$ with pre-compact range;
\item $L^\infty (\R,\X)$ denotes the space of all measurable
functions on $\R$ that are essentially bounded with usual norm $\|
g\| :={\rm ess}\sup_{t\in\R}\| g(t)\|$;
\item $\Gamma$ denotes the unit circle in the complex plane $\C$;
 \item If $B(t)$ is a bounded linear
operator on $\X$ for each $t\in\R$ that is strongly continuous in
$t$, then the operator of multiplication by $B(t)$ on $BC(J,\X)$,
denoted by ${\cal B}$, is defined by ${\cal B}g=Bg(\cdot )$ for each
$g\in BC(J,\X )$; \item The translation operator $S(\tau )$ is
defined to be $S(\tau )g (t)=g(t+\tau )$ for all $t\in \R$, $g\in
L^\infty (\R,\X)$; In particular, $S(1):= S$.
\end{enumerate}
\subsection{Function Spaces}
The biggest function space we consider in this paper is $L^\infty
(\R,\X)$ of all measurable functions that are essentially bounded
on $\R$ with ess-sup norm. We will identify each element in
$BC(\R,\X)$ with its equivalence class in $L^\infty (\R,\X)$, so
we may think of $BC(\R,\X)$ as a closed subspace of $L^\infty
(\R,\X)$.
\begin{definition} \rm
A function $f \in BC({\mathbb R}, {\mathbb X})$ is said to be {\it
almost automorphic} if for any sequence of real numbers $(s'_{n})$,
there exists a subsequence $(s_{n})$ such that
\begin{equation}\label{140}
\lim_{m \to \infty} \lim_{n \to \infty} f( t + s_{n} - s_{m}) = f(t)
\;
\end{equation}
for any $t \in {\mathbb R}$.
\end{definition}
The limit in (\ref{140}) means
\begin{equation}\label{41}
g(t) = \lim_{n \to \infty} f(t + s_{n})
\end{equation}
is well-defined for each $t \in {\mathbb R}$ and
\begin{eqnarray}\label{141}
f(t) = \lim_{n \to \infty} g(t - s_{n})
\end{eqnarray}
for each $t  \in {\mathbb R}$. The reader is referred to
\cite{ngu,ngu2,ngu3,minnaingu,bersie} and their references for more
information on this concept and results.

\begin{definition} \rm
A function $f \in BC({\mathbb R}, {\mathbb X})$ is said to be {\it
almost periodic} if for any sequence of real numbers $(s'_{n})$,
there exists a subsequence $(s_{n})$ and a function $g\in BC(\R,\X)$ such that
\begin{equation}\label{40}
\lim_{n \to \infty} f( t + s_{n}) = g(t)
\end{equation}
uniformly in $t\in\R$.
\end{definition}
It follows immediately from the definition that every almost
periodic function is uniformly continuous. The space of all almost
periodic functions on $\R$ taking values in $\X$ is denoted by
$AP(\X)$, so $AP(\X)\subset BUC(\R,\X)$. For more information on
almost periodic functions the reader is referred to
\cite{fin,levzhi}.

\begin{remark} \rm
Because of pointwise convergence the function $g$ is measurable but
not necessarily continuous. It is also clear from the definition
above that constant functions and almost periodic functions are
almost automorphic.
\end{remark}

If the limits in (\ref{41}) and (\ref{141}) are uniform on any compact subset $K
\subset {\mathbb R}$, we say that $f$ is compact almost automorphic.

\begin{theorem}
 Assume that $f$, $f_{1}$, and $f_{2}$ are almost
automorphic and $\lambda$ is any scalar, then the following hold
true.
\begin{enumerate}
\item $\lambda f$ and $f_{1} + f_{2}$ are almost automorphic,
\item $f_{\tau}(t) : = f(t + \tau)$, $t \in {\mathbb R}$ is almost
automorphic, \item $\bar{f} (t) : = f(- t)$, $t \in {\mathbb R}$ is
almost automorphic, \item The range $R_{f}$ of $f$ is precompact, so
$f$ is bounded.
\end{enumerate}
\end{theorem}

\begin{proof}
See \cite[Theorems 2.1.3 and 2.1.4]{ngu}, for proofs.
\end{proof}

\begin{theorem}
If $\{f_{n} \}$ is a sequence of almost automorphic ${\mathbb
X}$-valued functions such that $f_{n} \mapsto f$ uniformly on
${\mathbb R}$, then $f$ is almost automorphic.
\end{theorem}

\begin{proof}
see \cite[Theorem 2.1.10]{ngu}, for proof.
\end{proof}

\begin{remark} \rm
 If we equip $AA ({\mathbb X})$, the space of almost automorphic
functions with the sup norm
$$
\| f \|_{\infty} = \sup_{t \in {\mathbb R}} \| f(t) \|
$$
then it turns out to be a Banach space. If we denote $KAA({\mathbb
X})$, the space of compact almost automorphic ${\mathbb X}$-valued
functions, then we have
\begin{equation}
AP ({\mathbb X}) \subset KAA ({\mathbb X}) \subset AA ({\mathbb X})
\subset BC ({\mathbb R}, {\mathbb X}) \subset L^\infty (\R,\X).
\end{equation}
\end{remark}
\begin{theorem}
If $f \in AA({\mathbb X})$ and its derivative $f'$ exists and is
uniformly continuous on ${\mathbb R}$, then $f' \in AA({\mathbb
X})$.
\end{theorem}

\begin{proof}
See \cite[Theorem 2.4.1]{ngu} for a detailed proof.
\end{proof}

\begin{theorem}
Let us define $F: {\mathbb R} \mapsto {\mathbb X}$ by $F(t) =
\int_{0}^{t} f(s) ds$ where $f \in AA ( {\mathbb X})$. Then $F \in
AA({\mathbb X})$ iff $R_{F} = \{ F(t) | \ t \in {\mathbb R} \}$ is
precompact.
\end{theorem}
\begin{proof}
See \cite[Theorem 2.4.4]{ngu} for a detailed proof.
\end{proof}

\bigskip
As a big difference between almost periodic functions and almost
automorphic functions we remark that an almost automorphic
function is not necessarily uniformly continuous, as shown in the
following example due to B. M. Levitan:
\begin{example} The following function
$$f(t):=\sin \left( \frac{1}{2+\cos t +
\cos \sqrt{2}t} \right)$$ is almost automorphic, but not uniformly
continuous. Therefore, it is not almost periodic.
\end{example}

\section{A Spectral Theory of Functions}
Below we will introduce a transform of a function $g\in L^\infty
(\R,\X)$ on the real line that leads to a concept of spectrum of a
function. This spectrum coincides with the set of
$\overline{e^{isp(g)}}$ if in addition $g$ is  uniformly continuous.
Recall that $\Gamma$ denotes the unit circle in the complex plane.

\bigskip
Let $g\in L^\infty (\R,\X)$. Consider the complex function ${\cal
S}g(\lambda )$ in $\lambda \in \C \backslash \Gamma$ defined as
\begin{equation}\label{transform}
{\cal S}g(\lambda ):= R(\lambda ,S)g,\quad \lambda \in \C \backslash
\Gamma .
\end{equation}
Since $S$ is a translation, this transform is an analytic function
in $\lambda \in \C \backslash \Gamma $.
\begin{definition}\rm
The {\it circular spectrum} of $g\in L^\infty (\R,\X)$ is defined to
be the set of all $\xi_0\in \Gamma$ such that ${\cal S}g(\lambda )$
has no analytic extension into any neighborhood of $\xi_0$ in the
complex plane. This spectrum of $g$ is denoted by $\sigma (g)$ and
will be called for short {\it the spectrum of $g$} if this does not
cause any confusion. We will denote by $\rho (g)$ the set $\Gamma
\backslash \sigma (g)$.
\end{definition}
\begin{lemma} If $|\lambda | \not= 1$, then
\begin{equation}\label{norm of R(lambda ,S)}
\| R(\lambda , S) \| \le \frac{1}{|1-|\lambda||},
\end{equation}
and if $|\lambda |=1$, then $\lambda \in \sigma (S)$.
\end{lemma}
\begin{proof}
If $|\lambda|\not= 1$, since $\| S\| =1$, $\lambda \in \rho (S)$, so
we have
\begin{equation}
I=(\lambda -S)R(\lambda ,S)=\lambda R(\lambda ,S)-SR(\lambda ,S).
\end{equation}
Therefore, since $\| S\| =1$
\begin{eqnarray}
1&=& \| \lambda R(\lambda ,S)-SR(\lambda ,S)\| \nonumber \\
&\ge& | |\lambda | \cdot \| R(\lambda ,S) \| -\| S\| \cdot \|
R(\lambda ,S)\|| \nonumber \\
&=& | |\lambda | \cdot \| R(\lambda ,S) \| - \| R(\lambda ,S)\|| \nonumber\\
&=& | (|\lambda | -1) \cdot \| R(\lambda ,S) \| | .
\end{eqnarray}
This proves (\ref{norm of R(lambda ,S)}).

\medskip
If $|\lambda |=1$, then there is a real $\xi$ such that $\lambda
=e^{i\xi}.$ It can be easily seen that the function $g_\xi (t):=
e^{i\xi t}a$ for a fixed $0\not= a\in \X$, $t\in\R$ is an
eigenvector of $S$ in $L^\infty (\R,\X)$, so $\lambda \in \sigma
(S)$.
\end{proof}

Recall that if $A\in L(\X)$, then ${\cal A}$ denotes the operator
of multiplication by $A$, given by $(Ag)(t):=Ag(t), \forall t\in \R
$.
\begin{proposition}\label{pro 2.1}
Let $\{ g_n\}_{n=1}^\infty \subset L^\infty (\R,\X)$ such that
$g_n\to g\in L^\infty (\R,\X)$, and let $\Lambda$ be a closed subset
of the unit circle. Then the following assertions hold:
\begin{enumerate}
\item \ $\sigma (g)$ is closed. \item \ If $\sigma (g_{n}) \subset
\Lambda$ for all $n \in {\N}$, then $\sigma (g)\subset \Lambda $.
\item \ $\sigma ({\cal A}g)\subset \sigma (g)$ for all $A\in
L(\X)$. \item If $\sigma (g) =\emptyset$, then $g=0$.
\end{enumerate}
\end{proposition}
\begin{proof} i) The first assertion follows immediately from the definition.\\
 ii) The proof can be taken from that of \cite[Theorem 0.8, pp.21-22]{pru}.
In fact\footnote{We thank the referee for his several detailed
suggestions in the proof of part ii)}, assume that
$\lambda_0:=e^{i\theta_0} \in \Gamma \setminus \Lambda$. Since
$\Lambda$ is closed, we can choose $r>0$ such that, if $|\xi -
i\theta_0|<4r$, then ${\cal S}g_n(e^\xi)=R(e^{\xi},S)g_n$ is
analytic for all $n$. Let us choose a positive $\delta$ such that if
$|x |<\delta$, where $x,y$ are reals, and $\lambda =e^{x+iy}$, then
$$
\frac{1}{|1-|\lambda||} \le \frac{2}{|x|}.
$$
Notice that, if $0<|\Re \xi|<\delta$, where $\delta$ is the number
in the above,
\[
\|R(e^{\xi},S)g_n\|\le \|R(e^{\xi},S)\|\|g_n\| \le
\frac{1}{1-|e^{\xi}|}\|g_n\|\le \frac{2}{|\Re \xi|}\|g_n\|.
\]
We will show that, if $|\xi - i\theta_0|<r<\delta/4$,
\[
\|{\cal S}g_n(e^{\xi})\|\le \frac{16}{3r}\|g_n\|.
\]
Take the function $h(z)=(z-i\theta)(1+(z-i\theta)^2/4r^2)$. By
Cauchy's theorem, we have that, if $|\xi -i\theta_0|< r$,
\[
h(\xi)({\cal S}g_n(e^\xi)) =\frac{1}{2\pi i}\int
_{|z-i\theta_0|=2r}\frac{h(z){\cal S}g_n(e^z)}{z-\xi}dz,
\]
and that
\[
\|h(\xi)({\cal S}g_n(e^\xi))\|\le \frac{1}{2\pi} \int
_{|z-i\theta_0|=2r}2|\Re z|\frac{2\|g_n\|}{|\Re z|} \frac{|dz|}{|z
-\xi|}\le \frac{1}{2\pi}4\|g_n\|\frac{1}{r}2\pi 2r=8\|g_n\|.
\]
This implies that
\begin{eqnarray*}
\sup_{|\xi-i\theta_0|\le r}\|{\cal S}g_n(e^\xi)\|
=\sup_{|\xi-i\theta_0|= r}\|{\cal S}g_n(e^\xi)\| &\le& 8\|g_n\|
\sup_{|\xi-i\theta_0|= r}\left|\frac{1}{h(\xi)}\right| \\
&\le& 8\|g_n\|\frac{1}{(2r)3/4}\\ &=&\frac{16\|g_n\|}{3r}.
\end{eqnarray*}
Since $g_n \to g$, we can take $M$ such that $16\|g_n\|/3r\le M$ for
$n=1,2,\cdots,$. Set $U=\{\lambda=e^{\xi}:|\xi-i\theta_0|<r\}$. Then
$U$ is a neighborhood of $\lambda_0=e^{i\theta_0}$ and for $\lambda
\in U$
\begin{equation}\label{boundedness of Sgn}
\|{\cal S}g_n(\lambda)\|\le \frac{16\|g_n\|}{3r}\le M.
\end{equation}
If $|\lambda|\neq 1$, then
\[
\|{\cal S}g_n(\lambda)-{\cal S}g(\lambda)\|
=\|R(\lambda,S)(g_n-g)\|\le \frac{1}{|1-|\lambda||}\|g_n-g\|,
\]
so if $\lambda\in U\backslash \Gamma $
\begin{equation}\label{limit of Sg}
\lim_{n\to \infty}{\cal S}g_n(\lambda)={\cal
S}g(\lambda).\end{equation}

\medskip
By Vitali's Theorem \cite[Theorem A.5, p. 458]{arebathieneu}, the
uniform boundedness (\ref{boundedness of Sgn}) and (\ref{limit of
Sg}) yield that ${\cal S}g(\lambda)$ has an analytic extension to
$U$, that is, $\lambda_0\not\in \sigma (g)$.

\medskip\noindent
iii) The assertion is obvious.

\medskip\noindent
iv) If $\sigma (g)=\emptyset$, then ${\cal S}g(\lambda )$ is
analytic everywhere in $\C$. Moreover, by Lemma \ref{norm of
R(lambda ,S)} it should be bounded on $\C$. This shows that ${\cal
S}g(\lambda )$ is a constant. Again using Lemma \ref{norm of
R(lambda ,S)} we end up with this constant being zero. This yields
that $g=0$.
 \end{proof}

Below by ${\cal F}$ we denote one of the function spaces if not otherwise stated:
\begin{equation}\label{function spaces}
AP ({\mathbb X}) ; \  KAA ({\mathbb X}) ; \   AA ({\mathbb X}) ; \
KBUC(\R,\X); \ BC ({\mathbb R}, {\mathbb X}) ; \ L^\infty (\R,\X).
\end{equation}
\begin{lemma}\label{lem invariance of function spaces}
Let $\F$ be one of the function spaces $ AP ({\mathbb X}) ; \  KAA
({\mathbb X}) ; \   AA ({\mathbb X})$, and let ${\cal T}$ be a
bounded linear operator in $BC (\R,\X)$ that commutes with all
translations. Then, ${\cal T}$ leaves $\F$ invariant.
\end{lemma}
\begin{proof}
This lemma is obvious due to definitions of these function spaces
$\F$. Therefore, the detailed proofs are omitted.
\end{proof}

\begin{corollary}
Let $\Lambda$ be a closed subset of the unit circle. Then, the set
\begin{equation}\label{space Lambda (X)}
\Lambda_\F (\X) := \{ g\in {\cal F}|\ \sigma (g)\subset \Lambda \}
\end{equation}
is a closed subspace of $\F$.
\end{corollary}
\begin{proof}
It is easy to check that this is a linear subspace of $\F$.
Moreover, by (ii) of the above proposition, this space is closed.
\end{proof}

\begin{lemma}\label{lem spec of S on Lambda (X)}
Let $\Lambda$ be a closed subset of the unit circle. Then, the
translation operator $S$ leaves the space $\Lambda_\F (\X)$
invariant. Moreover,
\begin{equation}\label{spec of S on Lambda (X)}
\sigma (S|_{\Lambda_\F (\X)})=\Lambda .
\end{equation}
\end{lemma}
\begin{proof}
Since the function $g_\mu(t):=\mu^t$ is an eigenvector of
$S|_{\Lambda_\F (\X)}$ for each $\mu\in \Lambda$, it is clear that
$\sigma (S|_{\Lambda_\F (\X)})\supset \Lambda$. Now we prove the
converse. Let $\mu\in \Gamma$ but $\mu\not\in \Lambda$. We will show
that $\mu \in \rho (S|_{\Lambda_\F (\X)}).$ That is, for each $g\in
\Lambda _\F (\X)$ the equation
\begin{equation}\label{new 1}
\mu y-Sy=g,
\end{equation}
has a unique solution in $\Lambda_\F (\X).$

\medskip
Obviously, $R(\lambda ,S)g$ has an analytic extension in a
neighborhood of $\mu$. Moreover, note that $R(\lambda ,S)g\in
\Lambda_\F (\X)$ whenever $g\in\Lambda_\F (\X)$. Therefore,
(\ref{new 1}) has a solution $y_1:=\lim_{\lambda \to\mu } R(\lambda
,S)g\in\F$. This equation has a unique solution in $\F$. In fact,
suppose the homogeneous equation
$$
\mu y-Sy=0
$$
has a solution $y_0\in \Lambda_\F (\X)$. Then, since $\mu y_0=Sy_0$,
\begin{eqnarray*}
R(\lambda ,S)y_0 &=& \mu^{-1}R(\lambda ,S)Sy_0\\
&=&\mu^{-1} (\lambda R(\lambda ,S)y_0-y_0),
\end{eqnarray*}
so
\begin{eqnarray*}
(1 - \mu^{-1} \lambda )R(\lambda ,S)y_0&=& -\mu^{-1} y_0.
\end{eqnarray*}
Therefore,
\begin{eqnarray*}
 R(\lambda ,S)y_0&=& -\frac{\mu^{-1} y_0}{ 1 -
\mu^{-1} \lambda  }\\
&=& \frac{  y_0}{ \lambda - \mu  }.
\end{eqnarray*}
This shows that $\sigma (y_0)\subset \{ \mu\}$. And hence, $\sigma
(y_0)\subset \{ \mu\} \cap \Lambda=\emptyset$. By iv) of Proposition
\ref{pro 2.1} $y_0=0$, that is the uniqueness of the solution of the
homogeneous equation. This proves that $\mu \in \rho (S|_{\Lambda_\F
(\X)})$, and so the lemma is proved.
\end{proof}

Recall that for $u\in L^\infty (\R,\X)$, $sp(u)$ stands for the
Carleman spectrum, which consists of all $\xi \in \R$ such that the
Carleman transform of $u$, defined by
$$
\hat{u}(\lambda) := \begin{cases}
\begin{array}{ll}
\int_0^{\infty}e^{-\lambda t}u(t)dt     &(Re\lambda > 0)\\ \\
-\int^{\infty}_0 e^{\lambda t}u(-t)dt   &(Re\lambda < 0) ,
\end{array} \end{cases}
$$
has no holomorphic extension to any neighborhoods of $i\xi$. For
each $u\in BUC(\R,\X)$ we denote $\mathcal{M}_u :=
\overline{span\{S(\tau)u, \tau \in \R\}}\subset BUC(\R,\X)$. If
$u\in BUC(\R,\X)$, the Carleman spectrum of $u$ coincides with its
Arveson spectrum, defined by (see \cite[Lemma 4.6.8]{arebathieneu})
\begin{equation}\label{spec1}
i\ sp (u) = \sigma ({\mathcal D}_u ),
\end{equation}
where ${\mathcal D}_u$ is the infinitesimal generator of the
restriction of the group of translations $(S(t)|_{\mathcal{M}_u})_{t
\in \R}$ to the closed subspace $\mathcal{M}_u$.

\bigskip
The following lemma relates the spectrum $\sigma (g)$ with Carleman
spectrum of a uniformly continuous and bounded function.

\begin{lemma}
Let $g\in BUC(\R,\X)$. Then
\begin{equation}
\sigma (g)= \sigma (S|_{{\cal M}_g}).
\end{equation}
\end{lemma}
\begin{proof}
The definition of $\sigma (g)$ yields that
\begin{equation}\label{10}
\sigma (g) \subset \sigma (S|_{{\cal M}_g}).
\end{equation}
Now we prove the inverse inclusion. Let $\lambda_0\in \Gamma$ such
that $\lambda_0\not\in \sigma (g)$. By Lemma \ref{lem spec of S on
Lambda (X)} $\lambda_0 \in \rho (S|_{\Lambda_\F (\X)})$, where
$\Lambda := \sigma (g)$ and $\F :=BUC(\R,\X)$. From the definition
of ${\cal M}_g$ we can show that ${\cal M}_g\subset \Lambda_\F
(\X)$, and that
$$
R(\lambda_0 , S|_{\Lambda_\F (\X)}){\cal M}_g \subset {\cal M}_g.
$$
And thus $\lambda_0$ is in $\rho (S|_{{\cal M}_g})$. This proves the
inverse inclusion of (\ref{10}). The lemma is proved.
\end{proof}

\begin{corollary}\label{lem spec and WSM the}
Let $g\in BUC(\R,\X)$. Then
\begin{equation}
\sigma (g)= \overline{e^{i sp(g)}}.
\end{equation}
\end{corollary}
\begin{proof}
Since the translation group is bounded and strongly continuous in
$BUC(\R ,\X)$, by the Weak Spectral Mapping Theorem (see e.g.
\cite{engnag})
$$
\sigma (S|_{{\cal M}_g}) =\overline{e^{\sigma (\cD |_{{\cal
M}_g})}}.
$$
Therefore, the corollary follows immediately from the above lemma.
\end{proof}

\begin{remark} \rm
In general, for $g\in L^\infty (\R,\X)$ we do not know the relation
between the circular spectrum $\sigma (g)$ and its Carleman spectrum
$sp(g)$. We guess that $\sigma (g)$ may be larger than the set
$\overline{e^{isp(g)}}$.
\end{remark}

\section{Bounded Solutions of Difference Equations}
In this section we consider the existence of solutions in $\F$ as
one of the function spaces listed in (\ref{function spaces}) to
difference equations with continuous time of the form
\begin{equation}\label{1}
u(t)=B(t)u(t-1)+f(t),
\end{equation}
where $B(t)$ is a linear operator in a Banach space $\X$ that is
$1$-periodic, strongly continuous in $t$, and $f$ is in $\F$. We are
interested in finding conditions for the existence and uniqueness of
 solutions in $\F$ to (\ref{1}).

\bigskip
Below  we assume that $\F$ is one of the function spaces listed in (\ref{function spaces}). Then, under the above assumption on $B(t)$ the operator of multiplication by $B(t)$ in $L^\infty (\R,\X)$, denoted by
${\cal B}$, leaves $ \F  $ invariant. Therefore, it make sense to consider the restriction of $\cal B$ to $\F$, and to denote by $\sigma_\F ({\cal B})$ and $\rho_\F ({\cal B})$ the spectrum and resolvent set of this restriction, respectively. For simplicity we introduce e new notation:
\begin{equation}
\sigma_\F ({\cal B})\cap \Gamma =: \sigma_{\Gamma , \F}({\cal B}).
\end{equation}
When $\F$ is $L^\infty (\R,\X)$ we may use $ \sigma_{\Gamma  }({\cal B})$ instead of $ \sigma_{\Gamma , \F}({\cal B})$ if it does not cause any confusion.

\begin{lemma}
Let $f\in \F$, and let $u\in\F$ be a bounded solution of Eq. (\ref{1}).
Then, the following holds:
\begin{equation}\label{spectrum of bounded solution}
\sigma (u) \subset  \sigma_{\Gamma ,\F }({\cal B})  \cup \sigma (f).
\end{equation}
\end{lemma}
\begin{proof}
We consider the restrictions of $S$ and ${\cal B}$ to $\F$. First we prove the following identity for each $\lambda \not\in\Gamma$
\begin{equation}\label{2}
\lambda R(\lambda ,S)S(-1) =R(\lambda ,S) + S(-1).
\end{equation}
In fact, we have
\begin{eqnarray*}
R(\lambda ,S)(\lambda -S)S(-1) &=& S(-1)\\
R(\lambda ,S)(\lambda S(-1)-I) &=& S(-1)\\
\lambda R(\lambda ,S) S(-1) -R(\lambda ,S) &=& S(-1),
\end{eqnarray*}
so
\begin{eqnarray*}
\lambda R(\lambda ,S) S(-1) &=& R(\lambda ,S)+S(-1).
\end{eqnarray*}
Therefore, (\ref{2}) is valid.

\bigskip
Since $B(t)$ is $1$-periodic in $t$, we have $R(\lambda ,S) {\cal
B}={\cal B}R(\lambda ,S)$. As $u$ is a bounded solution, by
(\ref{2}) we have
\begin{eqnarray*}
R(\lambda ,S) u &=& R(\lambda ,S){\cal B}S(-1)u+R(\lambda ,S)f\\
&=& {\cal B} R(\lambda ,S)S(-1)u+R(\lambda ,S)f\\
&=& {\cal B} (\lambda^{-1}R(\lambda ,S) + \lambda^{-1}S(-1))u
+R(\lambda ,S)f,
\end{eqnarray*}
so, for each $\lambda\not\in\Gamma$,
\begin{eqnarray*}
\lambda R(\lambda ,S) u &=& {\cal B} (R(\lambda ,S) + S(-1))u +
\lambda R(\lambda ,S)f.
\end{eqnarray*}
Therefore, for each $\lambda\not\in\Gamma$,
\begin{eqnarray*}
\lambda R(\lambda ,S) u -{\cal B} R(\lambda ,S)u&=&   {\cal B}S(-1)u
+
\lambda R(\lambda ,S)f\\
(\lambda  -{\cal B}) R(\lambda ,S)u&=&  {\cal B} S(-1)u + \lambda
R(\lambda ,S)f.
\end{eqnarray*}
From this it follows that if $\lambda_0 \not\in \sigma_{\Gamma ,\F} ({\cal
B})$ and $\lambda_0\not\in \sigma (f)$, then $R(\lambda ,S)$ has an
analytic extension around a neighborhood of $\lambda_0$. This shows
(\ref{spectrum of bounded solution}).
\end{proof}

\begin{corollary}
Let $u$ and $f$ be in $\F$ that is one of the function spaces listed in (\ref{function
spaces}), and let $\sigma _{\Gamma,\F} ({\cal B})\cap \sigma (f)
=\emptyset$. Then, Eq.\ (\ref{1}) has no more than one solution
$u\in \Lambda_\F (\X)$, where $\Lambda :=\sigma (f)$.
\end{corollary}
\begin{proof}
It suffices to show that the homogeneous equation associated with
(\ref{1}) has no more than one solution in $\Lambda_\F (\X)$. Let
$u$ be such a solution of the homogeneous equation. By the above
lemma $\sigma (u) \subset \sigma _{\Gamma,\F} ({\cal B})$. Therefore,
$\sigma (u) \subset \sigma _{\Gamma,\F} ({\cal B})\cap \sigma
(f)=\emptyset ,$ so $u$ is the zero function.
\end{proof}

\begin{lemma}\label{lem spec of Qg}
Let $Q(t)$ be 1-periodic strongly continuous in $t$, $u\in L^\infty
(\R,\X)$, and let ${\cal Q}$ be the operator of multiplication by
$Q(t)$. Then
\begin{equation}\label{spec of Qg}
\sigma ({\cal Q}g) \subset \sigma (g).
\end{equation}
\end{lemma}
\begin{proof}
Since $\cal Q$ commutes with $R(\lambda ,S)$, we have
\begin{equation}
R(\lambda ,S){\cal Q}g ={\cal Q} R(\lambda ,S)g,
\end{equation}
so, $R(\lambda ,S){\cal Q}g $ has an analytic extension to a
neighborhood of $\lambda_0\in\Gamma$ whenever so does $R(\lambda
,S)g$.
\end{proof}
\begin{remark} \rm
By Lemma \ref{lem spec and WSM the}, (\ref{spec of Qg})  yields in particular that
if $g\in BUC(\R,\X)$, then
$$
e^{sp({\cal Q}g)}\subset \overline{e^{sp(g)}}.
$$
This spectral estimate was given in \cite{bathutrab}.
\end{remark}

\medskip
Recall that $\sigma_\F ({\cal B})$ denotes the spectrum of the restriction of ${\cal B}$ to $\F$, where $\F$ is one of the function spaces listed in (\ref{function
spaces}).
\begin{lemma}\label{spec of cal A on Lambda X}
Let $\F$ be one of the function spaces listed in (\ref{function
spaces}), and let $\Lambda$ be a closed subset of unit circle $\Gamma$. Then, under the above assumption on $B(t)$ the operator
${\cal B}$ leaves $\Lambda_\F (\X)$ invariant. Moreover,
\begin{equation}\label{spec of A on Lambda X}
\sigma ({\cal B}|_{\Lambda_\F (\X )}) \subset \sigma _\F ({\cal B}).
\end{equation}
\end{lemma}
\begin{proof}
For the first assertion we note that $\cal B$ leaves $\F$ invariant and commutes with $S$. Therefore, for $|\lambda |\not =1$, and $g\in \Lambda _\F (\X)$ we have
$$R(\lambda , S){\cal B}g ={\cal B}R(\lambda ,S)g,$$
so, $\sigma ({\cal B}g) \subset
\sigma ( g)\subset \Lambda $, that is, ${\cal B}g \in \Lambda_\F (\X)$.

\medskip
For the last
assertion suppose that $\lambda_0 \in \rho _\F ({\cal B})$. We will show
that $\lambda_0 \in \rho ({\cal B}|_{\Lambda _\F (\X )})$. In fact,
we have to show that for each $g\in \Lambda_\F (\X)$ the equation
\begin{equation}\label{4}
\lambda_0 u-{\cal B}u=g
\end{equation}
has a unique solution $u$ in $\Lambda_\F (\X)$. First, since
$\lambda_0\in \rho_\F ({\cal B})$, in $\F$ there exists a unique solution $R(\lambda_0 ,{\cal B})g=u$ of (\ref{4}). Therefore, if Eq. (\ref{4}) has a solution in $\Lambda _\F (\X)$, it cannot have more than one. Next,
since ${\cal B}$
commutes with $S$, we can prove that $R(\lambda _0, {\cal B})$
commutes with $R(\lambda ,S)$ for all $|\lambda|\not= 1$. And hence,
this yields that for all $|\lambda|\not= 1$,
\begin{eqnarray*}
R(\lambda ,S) u&=& R(\lambda ,S)R(\lambda_0 ,{\cal B})g\\
&=& R(\lambda_0 ,{\cal B})R(\lambda ,S)g.
\end{eqnarray*}
This shows that whenever $R(\lambda ,S)g$ has an analytic extension into a neighborhood of a complex number $\lambda_1\in\Gamma$, so does $R(\lambda ,S) u$. This yields that $\sigma (u) \subset \sigma (g)$, that is,  $u\in \Lambda _\F (\X)$. The lemma is proven.
\end{proof}
\begin{remark}\rm
As will be shown in the next section if $B(t)$ is good enough, e.g. the monodromy operator of a periodic evolution equation, the  spectrum and resolvent $\sigma_\F ({\cal B})$ and $\rho_\F ({\cal B})$ can be estimated independently of $\F$ (see Lemma \ref{lem 5.2} below).
\end{remark}

\begin{theorem}\label{the main for difference eq}
Let $f$ be in $\F$, where $\F$ is any of the function spaces listed
in (\ref{function spaces}). If
\begin{equation}
\sigma _\F ({\cal B}) \cap   \sigma (f) =\emptyset ,
\end{equation}
then, Eq. (\ref{1}) has a unique solution $u$ in $\F$ such that
\begin{equation}
\sigma (u) \subset \sigma (f).
\end{equation}
\end{theorem}
\begin{proof}
Consider the equation
\begin{equation}\label{5}
 u={\cal B}S(-1) u +f
\end{equation}
in $\Lambda_\F (\X)$, where $\Lambda :=\sigma (f)$. This equation is
equivalent to the following due to the commutativeness of ${\cal
B}|_{\Lambda_\F (\X)}$ and $S|_{\Lambda_\F (\X )}$
\begin{equation}\label{6}
( S|_{\Lambda_\F (\X)}-{\cal B}|_{\Lambda_\F (\X)})u =Sf.
\end{equation}
Moreover, the commutativeness of ${\cal B}|_{\Lambda_\F (\X)}$ and
$S|_{\Lambda (\X )}$ yields (see \cite{rud})
$$
\sigma ( S|_{\Lambda_\F (\X)})-{\cal B}|_{\Lambda_\F (\X)}) \subset
\sigma ( S|_{\Lambda_\F (\X)}) - \sigma ({\cal B}|_{\Lambda_\F
(\X)}) .
$$
By Lemmas \ref{lem spec of S on Lambda (X)} and \ref{spec of cal A on Lambda X}
$$
\sigma ( S|_{\Lambda_\F (\X)}) =\Lambda := \sigma (f), \quad \sigma ({\cal B}|_{\Lambda_\F
(\X)}) \subset \sigma_\F ({\cal B}).
$$
Therefore, by the theorem's assumption
$$
0\not\in \sigma (f) - \sigma_\F ({\cal B}) \supset
\sigma ( S|_{\Lambda_\F (\X)})-{\cal B}|_{\Lambda_\F (\X)})
$$
This shows that $0\in \rho (S|_{\Lambda_\F (\X)}-{\cal
B}|_{\Lambda_\F (\X)})$, that is, the operator $ \left( S|_{\Lambda_\F
(\X)}-{\cal B}|_{\Lambda_\F (\X)}\right)$ is invertible.
 In particular, this yields (\ref{6}), and thus (\ref{5}) has a
unique solution in $\Lambda_\F (\X)$. This proves the theorem.
\end{proof}

\section{Bounded mild solutions of periodic evolution equations}
As an application of the above result we consider the existence and
uniqueness of different classes of bounded mild solutions of
evolution equations of the form
\begin{equation}\label{eve}
\frac{du(t)}{dt}= A(t)u(t)+f(t), \quad t\in \R ,
\end{equation}
where $u(t) \in {\X}$, ${\X}$ is a complex Banach space, $A(t)$ is
a (unbounded) linear operator acting on ${\X}$ for every fixed $t
\in {\R}$ such that $A(t)=A(t+1)$ for all $t \in {\R}$,\linebreak
$f : {\R} \rightarrow {\X}$ is a bounded function. Under suitable
conditions the homogeneous equation associated with Eq.\
(\ref{eve}) is well-posed (see e.g.\ \cite{paz}), i.e., one can
associate with this equation an evolutionary process
$(U(t,s))_{t\ge s}$ which satisfies, among other things, the
conditions in the following definition.

\begin{definition}\rm
A family of bounded linear operators $(U(t,s))_{t\ge s}, (t,s\in
{\R})$ from a Banach space $\X$ to itself is called {\it 1-periodic
strongly continuous evolutionary process}
 if the following conditions are
satisfied:
\begin{enumerate}
\item $U(t,t)=I$ for all $t\in {\bf R}$,
\item $U(t,s)U(s,r)=U(t,r)$ for all $t\ge s\ge r$,
\item The map $(t,s)\mapsto U(t,s)x$ is continuous for every fixed $x\in {\X}$,
\item $U(t+1,s+1)=U(t,s)$ for all $t \ge s$ ,
\item $\| U(t,s)\| \le Ne^{\omega (t-s)} $ for some positive $N, \omega $ independent of
$t \geq s$ .
\end{enumerate}
\end{definition}

Recall that for a given 1-periodic evolutionary process
$(U(t,s))_{t\ge s}$ the following operator
\begin{equation}
P(t):= U(t,t-1), t \in {\bf R}
\end{equation}
is called {\it monodromy operator} (or sometime {\it period map,
Poincar\'e map}). Thus we have a family of monodromy operators. We
will denote $P:= P(0)$. The nonzero eigenvalues of $P(t)$  are
called {\it characteristic multipliers}. An important property of
monodromy operators is stated in the following lemma whose proof
can be found in \cite{hen,hinnaiminshi}.
\begin{lemma}\label{lem 5.2}
Under the notation as above the following assertions hold:
\begin{enumerate}
\item $P(t+1) = P(t)$ for all $t$; characteristic multipliers are
independent of time, i.e. the nonzero eigenvalues of $P(t)$ coincide
with those of $P$, \item $\sigma (P(t)) \backslash \{0\}=\sigma (P)
\backslash \{0\}$, i.e., it is independent of $t$, \item If $\lambda
\in \rho (P)$, then the resolvent $R(\lambda , P(t))$ is strongly
continuous, \item If ${\cal P}$ denotes the operator of
multiplication by $P(t)$ in any one of the function spaces $\F$ listed in
(\ref{function spaces}), then
\begin{equation}\label{5.3}
\sigma_\F ({\cal P}) \backslash \{0\} \subset \sigma (P)\backslash
\{0\}.
\end{equation}
\end{enumerate}
\end{lemma}
\begin{proof}
For (i)-(iii) the proofs can be found in \cite{hen,hinnaiminshi}.

\medskip
For (iv), for a fixed function space $\F$, note that by (i)-(iii), if $\lambda_0 \in \rho (P)$, then the operator of multiplication by $R(\lambda_0,P(t))$ leaves $\F$ invariant, so
$R(\lambda_0 , {\cal P})$ can be determined by $R(\lambda_0,P(t))$. Therefore, (\ref{5.3}) holds.
\end{proof}

We note that in the infinite dimensional case there does not
always exist a Floquet representation of the monodromy operator
$P$. And in general we do not know if by a "change of variables"
Eq.\ (\ref{eve}) can be reduced to an autonomous equation. In the
finite dimensional case, this can be done in the framework of the
Floquet Theory. If the Poincare map $P$ is compact, a partial
Floquet representation of $P$ may be used as in
\cite{hen,murnaimin2}. When $f$ is almost periodic it was
conjectured in \cite{vu} that the condition $\sigma _\Gamma
(P)\cap \sigma (f)$ is a sufficient condition for the existence
and uniqueness of an almost periodic mild solution $u$ to Eq.
(\ref{eve}) such that $\sigma (u)\subset \sigma (f)$. The
evolution semigroup method proposed in \cite{naimin} shows to be
working well to give a positive answer to the conjecture. For more
information about this we refer the reader to
\cite{naimin,bathutrab,murnaimin,hinnaiminshi}.

\bigskip
Recall that given a 1-periodic evolutionary process
$(U(t,s))_{t\ge s}$, the following formal semigroup associated
with it
\begin{equation}
(T^hu)(t):=U(t,t-h)u(t-h), \forall t\in {\R},
\end{equation}
where $u$ is an element of some function space, is called {\it
evolutionary semigroup} associated with the process $(U(t,s))_{t\ge
s}$. As is known, this evolution semigroup is strongly continuous at
each almost periodic function, or more generally at each bounded and
uniformly continuous function with pre-compact range. The strong
continuity of the evolution semigroup is essential in the evolution
semigroup method. However, it may not be strongly continuous at any
bounded and continuous function. Since an almost automorphic
function may not be uniformly continuous the extended conjecture of
Vu in \cite{vu} with almost automorphic $f$ is still open.

\bigskip
Below we will give a positive answer to the extended conjecture of
Vu with general $f\in BC(\R,\X)$ by applying the results in the
previous section.

Let $U(t,s)$ be a $1$-periodic strongly continuous evolutionary
process. We note that all results can be adjusted if the process
is $\tau$-periodic with any positive $\tau$. For each fixed
positive $h$ let us define an operator $G$ as follows
\begin{equation}\label{operator G}
Gg(t) := \int^{t}_{t-h} U(t,\xi )g(\xi )d\xi , \quad g\in L^\infty
(\R,\X), t\in\R .
\end{equation}
Note that this operator $G$ is well defined because of the strong
continuity of the process $(U(t,s))_{t\ge s}$. Moreover, $Gg\in
BC(\R,\X)$ for each $g\in L^\infty (\R,\X)$.

\begin{lemma}\label{lem spec of Gf}
Let $G$ be defined as above. Then the following assertions hold:
\begin{enumerate}
\item If $\F$ is one of the function spaces (\ref{function spaces}),
then $G$ leaves $\F$ invariant;
\item For each $g\in L^\infty (\R,\X)$,
\begin{equation}\label{spec of Gf}
 \sigma (Gg) \subset \sigma
(g).
\end{equation}
\end{enumerate}
\end{lemma}
\begin{proof}
By the $1$-periodicity of $(U(t,s))_{t,s\in\R}$, for all $t\in\R$
and all $f\in BC(\R,\X)$ we have
\begin{eqnarray*}
[SGg](t) &=& \int^{t+1}_{t-h+1} U(t+1,\xi )g(\xi )d\xi \\
&=& \int^{t}_{t-h} U(t+1,\eta +1 )g(\eta +1 )d\eta \\
&=& \int^{t}_{t-h} U(t,\eta )g(\eta +1)d\eta \\
&=& [GSg] (t),
\end{eqnarray*}
so $S$ commutes with $G$. This yields in particular that if $\F$ is
one of the function spaces $AP(\X),KAA(\X),AA(\X)$ the operator $G$
leaves $\F$ invariant. When $\F$ is one of the function spaces
$KBUC(\R,\X)$ or $BC(\R,\X)$, the invariance under $G$ can be easily
checked. Therefore, the first assertion follows.

\medskip\noindent
For the second assertion, the commutativeness of $G$ and $S$ yields
$$
R(\lambda ,S)Gg =GR(\lambda ,S)g,
$$
and $R(\lambda ,S)Gg$ has an analytic extension into a neighborhood
of $\lambda_0\in\Gamma$ whenever so does $R(\lambda ,S)g$. Finally,
this yields (\ref{spec of Gf}).
\end{proof}
Below we always assume that $\F$ is one of the function spaces in
(\ref{function spaces}). We consider the following semigroup
$(T^h_f)_{h\ge 0}$ of affine operators in $L^\infty (\R, \X)$: for each $h\ge 0 $ and $
f,g \in L^\infty (\R,\X),$
\begin{equation}
\left( T^h_f g\right)(t) := U(t,t-h)g(t-h) + \int^t_{t-h}U(t,\xi )
f(\xi )d\xi , \ \mbox{for almost all} \ t\in\R .
\end{equation}
Let $\Lambda$ be a closed subset of $\Gamma$, and let $f\in \Lambda_\F(\X)$. By Lemmas \ref{lem spec of Qg} and
\ref{lem spec of Gf} $T^h_f $ leaves $\Lambda_\F (\X)$ invariant.
Moreover, it forms a semigroup of operators in $\Lambda_\F (\X)$. In
fact, we will show that for any nonnegative $h,k$, $T^{h+k}_f
=T_f^hT^k_f$. To this end,
\begin{eqnarray*}
\left( T^{h+k}_f g\right)(t) &=& U(t,t-h-k)g(t-h-k) +
\int^t_{t-h-k}U(t,\xi ) f(\xi )d\xi \\
&=& U(t,t-h)\left ( U(t-h,t-h-k)g(t-h-k) +\int^{t-h}_{t-h-k}U(t,\xi
)g(\xi )d\xi  \right)\\
&&\hspace{1cm} +\int^t_{t-h}U(t,\xi )f(\xi )d\xi \\
&=& U(t,t-h) T^k_f g(t-h) +\int^t_{t-h}U(t,\xi )f(\xi )d\xi \\
&=& \left(T^h_f T^k_f g\right) (t).
\end{eqnarray*}
Recall that for a given $f\in BC(\R,\X)$, a function $u\in
BC(\R,\X)$ is a mild solution of a well-posed Eq.\ (\ref{eve}) if
\begin{equation}\label{mild sol}
u(t)=U(t,s)u(s)+\int^t_sU(t,\xi )f(\xi )d\xi , \ \mbox{for all } t\ge
s.
\end{equation}
\begin{lemma}
Let $f\in BC(\R,\X)$. Then, $u\in BC (\R,\X)$ is a mild solution of Eq.\ (\ref{eve}) if and only
if it is a common fixed point  for all operators of the semigroup
$(T^h_f)_{h\ge 0}$ in $BC(\R,\X)$.
\end{lemma}
\begin{proof}
This lemma is obvious.
\end{proof}

We are now ready to prove the main result of the paper:
\begin{theorem}\label{the main}
Let the homogeneous equation associated with Eq.\ (\ref{eve})
generate a $1$-periodic evolutionary process with monodromy operator
$P$, and let $f\in \F$, where $\F$ is one of the function spaces
$AP(\X), KAA(\X), AA(\X), KBUC(\R,\X), BC(\R,\X)$. Then, Eq.\
(\ref{eve}) has a unique mild solution $u$ in $\F$ such that
\begin{equation}
\sigma (u)\subset \sigma (f),
\end{equation}
provided that
\begin{equation}
\sigma (P) \cap \sigma (f)=\emptyset .
\end{equation}
\end{theorem}
\begin{proof}
By the above lemma, it suffices to show that the semigroup
$(T^h_f)_{h\ge 0}$ has a unique common fixed point in $\Lambda_\F
(\X)$, where $\Lambda =\sigma (f)$. By Theorem \ref{the main for
difference eq} the operator $T^1_f$ has a unique fixed point in
$\Lambda_\F (\X)$. We are going to show that it should be common for
all operators in the semigroup. In fact, let $u$ be the unique fixed
point for $T^1_f$. Then, for each $h\ge 0$,
\begin{eqnarray*}
T^1_fT^h_f  u&=& T^h_f T^1_fu\\
&=&T^h_f u
\end{eqnarray*}
so, $T^h_fu$ is another fixed point of $T^1_f$ such that $\sigma
(T^h_fu)\subset \sigma (f)$. By the uniqueness of the fixed point of
$T^1_f$ this yields that $T^h_fu=u$. And hence, $u$ is a common
fixed point (in $\Lambda_\F (\X)$) for the whole semigroup.
Therefore, $u$ is a mild solution of Eq. (\ref{eve}) such that
$\sigma (u)\subset \sigma (f)$. The uniqueness follows from the
uniqueness of the fixed point of $T^1_f$.
\end{proof}

\bigskip
Consider autonomous equations of the form
\begin{equation}\label{autonom eve}
\frac{du(t)}{dt}= Au(t)+f(t),\quad t\in\R ,
\end{equation}
where $A$ generates a $C_0$-semigroup $(T(t))_{t\ge 0}$, $f$ is an
$\X$-valued bounded and continuous function in $\F$ that is defined
in Theorem \ref{the main}. The case when $f$ is uniformly continuous
is well studied in \cite{pru,vu,vusch,hinnaiminshi}. Under the
assumption, the Poincare operator $P$ is nothing but $T(1)$. Note that in the autonomous case the operator of multiplication by a linear bounded operator $B$ in $BUC(\R,\X)$ leaves this space invariant, so in this case in addition to the function spaces listed in (\ref{function spaces}) we can add $BUC(\R,\X)$.
Therefore, the following corollary is valid.
\begin{corollary} Let $f$ be in $\F$ that is any function space in (\ref{function spaces}) or $BUC(\R,\X)$. Then,
Eq. (\ref{autonom eve}) has a unique mild solution $u\in \F$ such
that $\sigma (u) \subset \sigma (f)$ provided that
\begin{equation}
\sigma_\Gamma (T(1)) \cap \sigma (f) =\emptyset .
\end{equation}
\end{corollary}

\bigskip
Let us consider the perturbed equation (\ref{perturb peve}). We will
fix a closed subset $\Lambda \subset \Gamma$ and a function space
$\F$ being one of the function spaces
\begin{equation}\label{conti func spa}
AP(\X), KAA(\X),AA(\X), KBUC(\R,\X), BC(\R,\X).
\end{equation}
We assume that
\begin{enumerate}
\item[(H1)] $H(t,0)=0$, and $H(t,x)$ is $1$-periodic;
\item[(H2)]  There exists an
increasing function $l:\R^+ \to \R^+$ such that for each positive
$r$ and for all $x,y\in \{ \xi \in\X : \| \xi \| \le r\} $ and
$t\in\R$, the following holds
\begin{equation}\label{local Lip of H}
\| H(t,x)-H(t,y)\| \le l (r) \| x-y\| ;
\end{equation}
\item[(H3)]
The Nemytsky operator ${\cal H}$ acting in $\F$ induced by
$H$, that is, ${\cal H}g: t\mapsto H(t , g(t ))$ leaves
$\Lambda_\F(\X)$ invariant.
\end{enumerate}

\medskip
Before we proceed we recall an operator associated with the linear
equation (\ref{peve}). The operator $L$ associated with (\ref{peve})
is defined on $BC(\R,\X)$ with domain consisting of all $u\in
BC(\R,\X)$ such that there exists such a function $f\in BC(\R,\X)$
for which (\ref{mild sol}) holds. And in this case $Lu:=f$. As is
well known (see e.g. \cite{naimin,hinnaiminshi,liunguminvu}), $L$ is
a closed, single-valued operator acting on $BC(\R,\X)$.

\begin{lemma}\label{lem closedness of L}
Let $\F$ be one of the function spaces in (\ref{conti func spa}),
and let $\Lambda$ be a closed subset of the unit circle. Then the
restriction of the operator $L$ to $\Lambda_\F(\X)$ (denoted by
$L_{\F ,\Lambda}$) is closed.
\end{lemma}
\begin{proof}
Let $u_n\in D(L) \to u\in \Lambda_\F (\X)$ such that $Lu_n =f_n \to
f \in \Lambda_\F (\X)$. By definition, for each $n\in\N$,
\begin{equation}
u_n(t)=U(t,s)u_n(s)+\int^t_sU(t,\xi )f_n(\xi )d\xi , \ \mbox{for
all} t\ge s.
\end{equation}
For every fixed $(t,s)$ let $n$ tend to infinity, so we have
\begin{equation}
u(t)=U(t,s)u(s)+\int^t_sU(t,\xi )f(\xi )d\xi , \ \mbox{for all} t\ge
s.
\end{equation}
Since $u\in\Lambda_\F(\X)$ this shows that $u$ is in the domain of
the restriction of $L$ to $\Lambda_\F(\X)$. Therefore, the
restriction of $L$ to $\Lambda_\F(\X)$ is closed.
\end{proof}

In the sequel we will need the Inverse Function Theorem for
Lipschitz continuous mappings, that is the following lemma that can
be found as a slight modification of a well known result in \cite{mar,nit}.
\begin{lemma}\label{cor lip inverse function}
Let $T$ be a bounded operator from a Banach space $\X$ onto another
Banach space $\Y$ such that $T^{-1}$ exists as a bounded operator,
and let $\varphi :\X \to \Y$ is a Lipschitz continuous operator with
$$
Lip (\varphi ) < \| T^{-1}\| ^{-1}.
$$ Then,
$(T+\varphi )$ is invertible with a Lipschitz continuous inverse,
and
\begin{eqnarray*}
Lip((T+\varphi )^{-1}  ) \frac{1 }{\| T^{-1}\|^{-1}- Lip (\varphi
)};
\end{eqnarray*}
\end{lemma}

Below we assume that $\F$ is one of the function spaces in
(\ref{conti func spa}).
\begin{theorem}
Let the homogeneous equation associated with Eq.\ (\ref{eve})
generate a $1$-periodic evolutionary process with monodromy operator
$P$, $\Lambda$ be a closed subset of $\Gamma$ such that
\begin{equation}
\sigma (P) \cap \Lambda =\emptyset ,
\end{equation}
$\F$ be a any fixed space from (\ref{conti func spa}), and let $f\in \Lambda_\F (\X)$. Assume further that $H$ in
(\ref{perturb peve}) satisfies all conditions (H1), (H2), (H3).
Then, there exists a positive constant $\epsilon_0$ such that if
$\epsilon <\epsilon_0$, the perturbed equation (\ref{perturb peve})
has a bounded mild solution in $\Lambda_\F (\X)$ that is locally
unique.
\end{theorem}
\begin{proof}
We will use Lemma \ref{cor lip inverse function}. Let
$$M:=   \sup_{t\in\R}\|f (t)\| .$$
As shown in Lemma \ref{lem closedness of L} and in Theorem \ref{the
main}, the restriction of $L$ to $\Lambda_\F(\X)$ is closed and
invertible. Therefore, if we equip $\X_1:= D(L_{\F ,\Lambda })$ with
its graph norm, then $L^{-1}_{\F,\Lambda}: \X_2:= \Lambda_\F (\X)
\to \X_1:= D(L_{\F ,\Lambda })$ is bounded. Let us denote
$$
\rho := \|L^{-1}_{\F,\Lambda}\| .
$$
We define the cut-off mapping
$$
{\cal H}_M (\phi )=
\begin{cases}
{\cal H}( \phi ), \ \ \forall \phi   \ \mbox{with} \ \| \phi \| \le 2\rho M \\
{\cal H}( \frac{2\rho M}{\| \phi \| }\phi ), \ \ \forall \phi  \
\mbox{with} \ \| \phi \| >2\rho M .
\end{cases}
$$
As is shown in \cite[Proposition 3.10, p.95]{web1}, in $B(2\rho M
):=\{ \phi \in\X_1 : \| \phi \| \le 2\rho M \}$, ${\cal H}_M(\cdot )$ is
globally Lipschitz continuous (in the new graph norm) with
$${Lip} ({\cal H}_M ) \le 2 {Lip}({\cal H}|_{ B(2\rho M ) }),
$$
where ${Lip} (R)$ denotes the Lipschitz coefficient of an operator $R:\X_1 \to \X_2$, so ${\cal H}_M$ satisfies
$$
Lip( {\cal H}_M ) \le 2l(2\rho M ).
$$
Since
\begin{equation}\label{112}
\lim_{t \downarrow 0}\frac{\rho  }{1-\rho t }=\rho >0
\end{equation}
we can choose $\epsilon_1$ so that
\begin{equation}
 \frac{\rho   }{1-\rho\epsilon 2l(2\rho M) } <2\rho
\end{equation}
for all $\epsilon <\epsilon_1$. By Lemma \ref{cor lip inverse
function}, if we choose $\epsilon_2$ such that
\begin{equation}
\epsilon_2 =\frac{\rho}{2l(2\rho M)},
\end{equation}
then, since $Lip(\epsilon {\cal H}_M)<\rho =
\|L^{-1}_{\F,\Lambda}\|$, the operator $L_{\F,\Lambda} + \epsilon
{\cal H}_M$ is invertible for all $\epsilon <\epsilon_2$. Finally,
if we choose $\epsilon_0= \min (\epsilon_1,\epsilon_2)$, then
$(L_{\F,\Lambda} + \epsilon {\cal H}_M)^{-1}$ exists and (\ref{112})
holds. Note that ${\cal H}_M (0)=0$. Therefore, if we let
$T:=L_{\F,\Lambda}$, $\varphi =\epsilon {\cal H}_M$ for $\epsilon
<\epsilon_0$,  then, by the above corollary
\begin{eqnarray*}
\| (L_{\F,\Lambda} - \epsilon {\cal H}_M)^{-1} f\| &= & \| (T+\varphi )^{-1} f\| \\
&=&  \| (T+\varphi )^{-1} f - (T+\varphi
)^{-1} (0)\| \\
 &\le& \frac{M }{\| T^{-1}\|^{-1}- Lip (\varphi
)}\\
&=&
\frac{M }{\rho^{-1}- \epsilon _0 2l(2\rho M))}\\
&=& \frac{M \rho }{1- \rho \epsilon _0 2l(2\rho M))}\\
&\le& 2\rho M.
\end{eqnarray*}
This shows that if $w:=(L_{\F,\Lambda} - \epsilon {\cal H}_M)^{-1}
f$, then $ (L_{\F,\Lambda} - \epsilon {\cal H}_M)w=f$, and $\| w\|
\le 2\rho M$. By the definition of ${\cal H}_M$, if $\| w\| \le
2\rho M$, then ${\cal H}_Mw={\cal H}w$. Finally, this yields that
$$
(L_{\F,\Lambda} - \epsilon {\cal H})w=f,
$$
that is, $w$ is a mild solution of (\ref{perturb peve}). The theorem
is proved.
\end{proof}

\section{Concluding remarks}
The choice of $1$-periodicity for the evolution equations in this
paper does not restrict the generality of the obtained results.
However, when dealing with $\tau$-periodic evolution equations the
concept of circular spectrum should be adjusted. Instead of using
the transform $R(\lambda,S)$ we use $R(\lambda , S(\tau ))$. If we
denote this spectrum by $\sigma^\tau (g)$, then the relation between
the Carleman spectrum and this circular spectrum can be established
in the following for $g\in BUC(\R,\X)$ via the Weak Spectral Mapping
Theorem
$$
\sigma^\tau (g) =\overline{e^{i\tau sp (g)}}.
$$
We may extend a little the statements of the results by refining
the classes of functions for $\F$ to be taken as in
\cite{bathutrab,hinnaiminshi,murnaimin,naimin,vusch}.

\end{document}